\documentclass{amsart}

\def\uh{\hat{u}}   

\def\eg{{\it e.g.\ }} \def\ie{{\it i.e.\ }}
\def\Real{{\mathbb R}}
\def\tQ{{\tilde Q}}
\def\EBj{E_{\mathcal B(j)}}
\def\EBjp1{E_{\mathcal B(j+1)}}
\def\Ai{\mbox{Ai}}

\theoremstyle{plain}
\newtheorem{thm}{Theorem}
\newtheorem{lem}{Lemma}

\theoremstyle{remark}

\theoremstyle{definition}

\begin{document}
\pagestyle{myheadings}
\markright{Fabian Waleffe \hfill  September 14, 2004 \hfill}

\title{On some dyadic models of the Euler equations}
\author{Fabian Waleffe}
\address{Departments of Mathematics and Engineering Physics,
University of Wisconsin, Madison, WI 53706}
\address{CNLS, CCS-2 and IGPP, Los Alamos National Laboratory,
Los Alamos, NM 87544. [LA-UR-04-7222]}
\email{waleffe@math.wisc.edu}
\date{\today}

\begin{abstract}
Katz and Pavlovic recently proposed a dyadic model of the Euler 
equations for which they proved finite time blow-up in the 
$H^{3/2+\epsilon}$ Sobolev norm. It is shown that their model 
can be reduced to the dyadic inviscid Burgers equation where nonlinear 
interactions are restricted to dyadic wavenumbers. The inviscid
Burgers equation  exhibits finite time blow-up in
 $H^{\alpha}$, for $\alpha \ge 1/2$,
  but its dyadic restriction is even more singular, exhibiting
blow-up for any $\alpha >0$.
 Friedlander and Pavlovic developed a closely related model
for which they also prove finite time blow-up in
$H^{3/2+\epsilon}$. Some inconsistent assumptions in the 
construction of their model are outlined.
 Finite time blow-up in the $H^{\alpha}$ norm, with $\alpha >0$,
is proven for a class of models that includes all those models.
An alternative shell model of the 
Navier-Stokes equations is discussed.
\end{abstract}

\keywords{Euler equations, Burgers equation, Navier-Stokes equations, finite time blow-up}
\subjclass{}
\maketitle

\section{Introduction}

The Euler equations of compressible gas dynamics typically develop
singularities in finite time from smooth initial conditions.
The singularities are known as shocks in that context.
The simplest example of this phenomenon
is provided by the {\it inviscid Burgers,} or {\it traffic flow},
equation
\begin{equation}
\frac{\partial u}{\partial t} + u \frac{\partial u}{\partial x} =0,
\label{burgers}
\end{equation}
where $u=u(x,t)$ and $x$ are real, $t >  0$.
 Let $\zeta(x,t)=\partial u/\partial x$, then from (\ref{burgers})
\begin{equation}
\frac{\partial \zeta}{\partial t} + u \frac{\partial \zeta}{\partial x}
= - \zeta^2.
\label{bslope}
\end{equation}
By the method of characteristics, we obtain
\begin{equation}
\frac{d \zeta}{dt} = - \zeta^2,
\label{blow}
\end{equation}
for  $\zeta=\zeta(X(t),t)$, 
along the characteristic lines $x=X(t)$ such that $dX/dt=u(X(t),t)$
with $du/dt=0$ for $u=u(X(t),t)$.
Equation
(\ref{blow}) exhibits blow-up in finite time for negative initial
conditions. Its solution is 
\begin{equation}
\zeta= \frac{\zeta_0}{1+\zeta_0 t}
\label{blow2}
\end{equation}
which is singular at $t=-1/\zeta_0 >0$ if $\zeta_0<0$. 

This singularity leads to finite time blow-up in the $H^{\alpha}$
Sobolev norm for $\alpha \ge 1/2$. 
The $H^{\alpha}$ norm of $u$, denoted $\|u\|_{H^{\alpha}}$, 
is defined in terms of the Fourier transform $\uh(k,t)$ 
of $u(x,t)$ as 
\begin{equation}
\|u\|^2_{H^{\alpha}} = 
\int_{\Real} \left(1 + |k|^{2 \alpha} \right) |\uh(k,t)|^2 dk.
\label{Halpha}
\end{equation}
Now, the method of characteristics provides the solution of
(\ref{burgers}) with $u(x,0)=f(x)$ in
implicit parametric form as $u=f(\eta)$, $x=\eta + t f(\eta)$,  thus
\begin{equation}
\begin{split}
2 \pi \uh(k,t) &=  \int_{\Real} u(x,t) e^{-ikx} dx\\
&= \int_{\Real} f(\eta) \, 
e^{-ik \left[\eta + t f(\eta)\right]} \, \left[1 + t f'(\eta)\right] d\eta\\
&= \frac{1}{ik} \int_{\Real} f'(\eta) \, 
e^{-ik \left[\eta + t f(\eta)\right]} \, d\eta,
\end{split}
\label{uhateta}
\end{equation}
after integration by parts, where $f'(\eta)=df/d\eta$. 
The singularity occurs when $1+ t f'(\eta_0)=0$ for some $\eta_0$.
For simplicity, assume that there is only one such $\eta_0$
with $f'(\eta_0)<0$,  $f''(\eta_0)=0$ and $f'''(\eta_0) =f'''_0 \neq 0$,
corresponding to the most negative initial slope. 
The asymptotic behavior of $\uh(k,t)$ as
$k \to \infty$, at the time of singularity, then follows by the method
of stationary phase as 
\begin{equation}
\begin{split}
\uh(k,t) &\sim \frac{1}{2 \pi i k} \, f'(\eta_0) \,
e^{-ik \left[\eta_0 + t f(\eta_0)\right]} 
\int_{\Real} e^{-ik t f'''_0 s^3/3!} \, ds \\
& = \frac{\Ai(0) }{ik} \left(\frac{2}{kt f'''_0}\right)^{1/3} 
\, f'(\eta_0) \, e^{-ik \left[\eta_0 + t f(\eta_0)\right]} , 
\end{split}
\label{uhSP}
\end{equation}
where $\Ai(0) = (2\pi)^{-1} \int_{\Real} e^{i y^3/3}dy$ is the
Airy function of the first kind evaluated at the origin. 
For example, if $f(\eta)=-\eta e^{-\eta^2}$  then 
$\uh(k,1) \sim i \, 3^{-1/3} \Ai(0) \,k^{-4/3}$. 
Since $|\uh(k,t)| = O(k^{-4/3})$ as $k \to \infty$, it follows
that the $H^{\alpha}$ norm (\ref{Halpha}) will diverge if 
$\alpha \ge 5/6$. For non-generic initial conditions
such that $f'''_0=0$, finite time blow-up can occur for smaller
$\alpha$. Smooth initial conditions where $f(\eta)=-\eta$
 in a finite neighborhood of the origin, for instance, lead to  
finite time blow-up for $\alpha \ge 1/2$, since $|\uh(k,1)|=O(k^{-1})$
as $k \to \infty$  in such cases.

It is an open question 
whether singularities can develop in finite-time for the Euler
equations when the flow is incompressible. This fundamental question
is related to the phenomenon of turbulence in fluid flows
since classical turbulence phenomenology, as well as experimental
data,  suggest that the energy dissipation rate in turbulent flows 
tends to a strictly positive constant as the fluid viscosity tends 
to zero \cite{Frisch}.

 Katz and Pavlovic \cite{KP04}
and Friedlander and Pavlovic \cite{FP04} recently proposed dyadic models of
the incompressible Euler equations and proved that their models exhibit
finite-time blow-up in the $H^{3/2+\epsilon}$ norm.
Their proof is restricted to $0 < \epsilon <2/3$.
We show that the Katz-Pavlovic model is in fact equivalent to a dyadic model
of the inviscid Burgers equation. The latter exhibits finite-time
blow-up in $H^{\alpha}$, for all $\alpha >0$,
 much stronger even than the unrestricted
inviscid Burgers equation which shows blow-up 
for $\alpha \ge 1/2$. Some
inconsistencies in the development of 
the Friedlander and Pavlovic model are outlined.
It is shown that their model and the dyadic
inviscid Burgers model belong to a class of models
that exhibit finite-time blow-up in $H^{\alpha}$, for all $\alpha >0$, 
even for weak nonlinear couplings. 

\section{Dyadic Inviscid Burgers equation}

After the rescaling $u \to 2 u$,
the inviscid Burgers equation (\ref{burgers}) can be written in 
conservative form as
\begin{equation}
\frac{\partial u}{\partial t} + \frac{\partial (u^2)}{\partial x} =0.
\label{burg2}
\end{equation}
Consider periodic solutions in  $-2^{-j_0}\pi \le x \le 2^{-j_0}\pi$, for
 some $j_0  \in \mathbb{Z}$.
 Let $\uh_l(t)$ be the $l$-th Fourier coefficient of 
\begin{equation}
u(x,t) = \sum_{l=-\infty}^{\infty} \uh_l(t) \, e^{i k_l x} 
\end{equation}
where the wavenumber $k_l= 2^{j_0} l$. 
The Fourier transform of equation (\ref{burg2}) reads
\begin{equation}
\frac{d \uh_l}{d t} 
= - i k_l \sum_{n=-\infty}^{\infty} \uh_n \uh_{l-n},
\label{FTburg}
\end{equation}
where the time dependency of $\uh_l=\uh_l(t)$ has been kept
implicit. Consider odd initial conditions so the solution remains
odd for all times, $u(x,t)=-u(-x,t)$, $\forall t \ge 0$. Since
$u(x,t)$ is real, this anti-symmetry implies that $\uh_l(t)$ is pure imaginary,
hence, let $\uh_l(t)= i v_l(t)$ where $v_l(t)$ is real and odd 
in $l$, so $v_l(t)=-v_{-l}(t)$, and we need only consider $l \ge 1$.
 Equation (\ref{FTburg}) becomes
\begin{equation}
\frac{d v_l}{dt} 
=  k_l \sum_{n=-\infty}^{\infty} v_n v_{l-n}.
\label{FTburg2}
\end{equation}
Now, artificially restrict $v_l$ to $l= 2^m$, with $m \ge 0$, corresponding
to wavenumber $k_l=2^{j_0} l =  2^{j_0+m}$. Let $j=j_0+m$
and define $a_j(t)=v_{l}(t)$. Using $v_n=-v_{-n}$, equation (\ref{FTburg2}) becomes
\begin{equation} \begin{split}
\frac{d a_j}{d t} 
=  2^j \Bigl( a_{j-1}^2 - 2 a_j a_{j+1} \Bigr), \quad j>j_0 \\
\frac{d a_{j_0}}{d t} =  -2^{j_0+1} a_{j_0} a_{j_0+1}. \qquad \qquad
\end{split}
\label{dyburg}
\end{equation}

\section{Reduction of the Katz-Pavlovic model}
\label{KPredux}

The Katz-Pavlovic (KP) model is based, formally, on a wavelet expansion of a 
scalar function $u(x,t)$, with $x$ $\in$ $\Real^3$,  over the set of 
dyadic cubes in $\Real^3$.
 This is the set of all cubes having
sidelength $2^{-j}$ with corners  on the lattice $2^{-j} \mathbb{Z}^3$. 
If $Q$ is a dyadic cube of sides $2^{-j}$, then $\tQ$ is the unique
parent cube of sidelength $2^{-j+1}$ containing $Q$ and $\mathcal{C}^1(Q)$ 
is the set of all $2^3$ children of $Q$, each having sidelength
$2^{-j-1}$. For $m \ge 1$, $\mathcal{C}^m(Q)$ denotes the set
of  $m$th order grandchildren of $Q$, \ie the set of all cubes of 
sidelength $2^{-j-m}$ that are contained in $Q$.

The KP equations are not derived from the Euler equations,
they are chosen to mimic the energy conserving quadratic nonlinearity 
of the Euler equations. The nonlinear interactions are restricted to 
local interactions in wavelet space and designed to push energy to smaller scales.
The KP model equations for the amplitude $u_Q(t)$ of the wavelet
localized at cube $Q$ of sidelength $2^{-j}$ is \cite[eqn.\ (2.5)]{KP04}
\begin{equation}
\frac{du_Q}{dt} =  2^{5j/2} u^2_{\tQ} - 
2^{5(j+1)/2} u_Q \sum_{Q'  \in \mathcal{C}^1(Q)} u_{Q'}.
\end{equation}
The factors $2^{5j/2}$ were chosen based on the scaling properties
of the Euler nonlinearity in $\Real^3$ (see \cite{KP04} for details).

The KP model is a tree, with each mode having 8 children, however it
is a simple tree where each edge at the same level has the same
weight and each cube interacts only with its unique parent and its own
8 children.
Therefore if, for any cube $\tQ_0$ of sidelength  $2^{-j_0}$,
the initial conditions are such that all $m$th order grandchidren 
of $\tQ_0$ have the same amplitude, \ie 
if $u_Q(t_0) = u_j(t_0)$ for all cubes $Q$ of sidelength $2^{-j}$
in $\mathcal{C}^m(\tQ_0)$ (so $j=j_0+m$), for all $m>1$,
 then they remain so for all $t > t_0$. 
For such initial conditions, the dynamics of the $\tQ_0$ branch
reduces to a chain model
\begin{equation}
\frac{du_j(t)}{dt} =  2^{5j/2} u^2_{j-1} - 
2^3 2^{5(j+1)/2} u_j u_{j+1},  
\label{KPchain}
\end{equation}
for all $j>j_0$.
The factor $2^3$ in the 2nd term on the right-hand side arises
from the fact that each dyadic cube has $2^3$ children 
$Q'$ in $\mathcal{C}^1(Q)$ and we now have $u_{Q'}(t)= u_{j+1}(t)$ for 
all dyadic cube $Q'$ with sidelength $2^{-j-1}$.
Now if the initial conditions are also zero for all cubes 
with $j \le j_0$ except for the cube $\tQ_0$, then the complete
dynamics is described by the chain (\ref{KPchain}) for $j>j_0$
together with 
\begin{equation}
\frac{du_{j_0}}{dt} =   - 2^3 2^{5(j_0+1)/2} u_{j_0} u_{j_0+1},
\label{j0mode}
\end{equation}
for the amplitude of the root cube $\tQ_0$.

For such initial conditions, the total energy reads
\begin{equation}
E= \sum_{m=0}^{\infty} 2^{3m} u^2_{j0+m} = 
2^{-3 j_0}\sum_{j=j_0}^{\infty} 2^{3j} u^2_j.
\label{Etot}
\end{equation}
The total energy is conserved by the dynamics (\ref{KPchain}),
(\ref{j0mode}). Our interest is in finite energy solutions.

Define $a_j = 2^{3j/2} u_j$ then the energy reads 
$E = 2^{-3 j_0} \sum_{j=j_0}^{\infty} a_j^2$ and the equations 
(\ref{KPchain}), (\ref{j0mode}), become
\begin{equation}
\begin{split}
\frac{1}{8}\frac{d a_j}{d t} =  2^j  a_{j-1}^2 - 2^{j+1} a_j a_{j+1},
\quad j >j_0, \\
\frac{1}{8} \frac{d a_{j_0}}{d t} =  -2^{j_0+1} a_{j_0} a_{j_0+1}. \qquad \qquad
\end{split}
\end{equation}
These equations are identical to the dyadic inviscid Burgers 
equations (\ref{dyburg}), after a rescaling of
time to eliminate the $1/8$ factor on the left-hand side.

\section{Friedlander-Pavlovic Model}

The Friedlander-Pavlovic model consists of the chain of ODEs
\cite[eqn.\ (3.10)]{FP04}
\begin{equation}
2\frac{d a_j}{d t} 
=  2^{5j/2}  a_{j-1}^2 - 2^{5(j+1)/2} a_j a_{j+1}.
\label{FPmod}
\end{equation}
This scalar chain model is deduced from a more complex vector model.
The latter vector model is constructed by drawing elements from the
Katz-Pavlovic model, the Dinaburg-Sinai model \cite{DS01} and the general
philosophy of shell models designed to study some features of
homogeneous turbulent flows (see \cite{Bif02} for a recent review of shell models).

The Dinaburg-Sinai model has been discussed elsewhere
\cite{Wds04}. Briefly, it is a model of the incompressible 
Navier-Stokes equations in an unbounded domain, 
deduced by assuming that the nonlinearity is dominated by highly
non-local interactions in Fourier space. In that model, small scales
do not interact directly with each other, they evolve due to the
distortion by an infinitely large scale background flow with uniform
gradient. Strictly speaking, this leads to a linear model of the 
Navier-Stokes equations since the large scale flow would have to be
determined externally. However, Dinaburg and Sinai inconsistently
define the infinitely large scale gradient as the linear superposition
of the small scale gradients at the origin.
This leads them to a quasilinear system of equations for the small
scales (see \cite{Wds04} for details). The derivation of the
model is physically and mathematically inconsistent, leading
in particular to a lack of conservation of energy since the infinitely 
large scale flow contains an infinite amount of energy,
but it is possible to  prove that the resulting system of 
equations has smooth solutions for limited classes of initial
conditions \cite{DS01}, \cite{DS03}.

Friedlander and Pavlovic begin with the shell model approach of 
considering an infinite sequence of nested shells $S_j$ in Fourier
space with exponential spacing:
$S_j=\{k \in \Real^3: 2^{j-1} \le |k| < 2^{j+1}\}$,
where  $j$ is an integer. 
They consider a single time-dependent wavevector $k_j(t)$ $\in \Real^3$
 with associated  velocity $v_j(t)$ $\in \Real^3$ 
in each shell $S_j$ and define the nonlinear 
interactions following the Katz-Pavlovic local interactions model
but with some of the structure of the Dinaburg-Sinai non-local interactions 
model. As a result, their wavevector evolution 
corresponds to the distortion by an infinitely large scale flow 
with uniform gradient, but the latter is now defined as the symmetric 
part of the gradient due to the next smaller scale at the origin
(\cite[eqn.\ (2.35)]{KP04}). This is still physically 
and mathematically inconsistent from a modeling and asymptotics
 standpoint. It is unclear here why they choose to symmetrize 
their matrix $B^T$ which should represent the ``large scale'' gradient.
However, since they assume that each wavevector is only distorted by
the next higher wavenumber, the wavevector dynamics 
is relatively simple. 
This allows them to construct special classes of initial conditions for 
which the wavevector dynamics and the entire vector structure of
the model disappear, except for a  dynamically
insignificant factor of 2 on the left-hand side of (\ref{FPmod}).
The symmetrization of $B^T$ appears key to that scalar reduction, since
if it was defined as a velocity gradient,   
$B^T=\kappa^{(j+1)} \left(v^{(j+1)}\right)^T$,
then  eqn.\ (3.3) in \cite{KP04} would equal $\alpha^{(j+1)}/\sqrt{2}$ instead
of  $\alpha^{(j)}/{2}$.

The Friedlander-Pavlovic scalar model (\ref{FPmod}) is structurally similar to 
the  dyadic inviscid Burgers model (\ref{dyburg}), 
which we showed to be equivalent to the
Katz-Pavlovic model for a special class of initial
conditions. However, the interactions coefficients in (\ref{FPmod}) 
are $2^{5j/2}$ instead of $2^j$. This appears to be another
inconsistency of the Friedlander-Pavlovic model since they consider
a single wavevector of magnitude $2^j$ in each shell $S_j$ 
but their estimate for the nonlinearity,
which is valid for wavelets in $\Real^3$, is an 
overestimate for a single Fourier mode.

\section{Finite-time blow-up in $H^{\alpha}$}

Here, we consider the class of chain models
\begin{equation}
\begin{split}
\frac{d a_j}{d t} 
= & \lambda^j  a_{j-1}^2 - \lambda^{j+1} a_j a_{j+1}, \quad j > j_0 \\
\frac{d a_{j_0}}{d t} 
= & \qquad - \lambda^{j_0+1} a_{j_0} a_{j_0+1}, \\
\end{split}
\label{cmod}
\end{equation}
for some integer $j_0$, where $\lambda > 1$ and the $a_j$'s are real.
This includes the dyadic Burgers model (\ref{dyburg})
when $\lambda=2$  and the  Friedlander-Pavlovic model (\ref{FPmod}) when
 $\lambda=2^{5/2}$.  

The total energy of system (\ref{cmod}) is defined as 
\begin{equation}
E = \sum_{j=j_0}^{\infty} E_j < \infty, \qquad E_j=a_j^2.
\end{equation}
It is easily verified that the total energy
is preserved by the dynamics (\ref{cmod}), thus  $E(t)=E(0)$,
for all $t>0$.
It is easily verified also that if the initial conditions are
non-negative, $a_j(0)\ge 0$, $\forall j$, then 
$a_j(t) \ge 0$, $\forall j$, $\forall t>0$, since  $da_j/dt \ge 0$
if $a_j(t)=0$. We consider such non-negative initial conditions
hereafter. 

A key characteristic of model (\ref{cmod}) with non-negative
initial data is that the energy flux is strictly toward higher
wavenumbers (\ie smaller scales). More precisely, let
\begin{equation}
\EBj= \sum_{l=j}^{\infty} E_l
\label{EBj}
\end{equation}
then from (\ref{cmod})
\begin{equation}
\frac{d}{dt}\EBj=  2 \lambda^j E_{j-1} a_j \ge 0,
\label{Eflux}
\end{equation}
since $a_j(t) \ge 0$.

The following proof of finite-time blow-up
is a generalization of Lemma 5.2 and Theorem 5.3 in
\cite{FP04} which are themselves a rewrite of lemma 3.2.2 and 
corollary 3.2.3 in \cite{KP04}.

\begin{lem}
For some $q$ with $\lambda^{-2} < q < 1$, and for some $j$ sufficiently large, 
if  
\begin{equation}
\EBj(t_0) \ge q^j
\label{qj}
\end{equation}
then 
there is $\rho$ with $\left(\lambda^2 q \right)^{-1/2} < \rho <1$, 
and $t$ in $[t_0, t_0+\rho^j]$,  such that 
$$\EBjp1(t) \ge q^{j+1}.$$ 
\label{lemma1}
\end{lem}

\begin{proof}
The proof proceeds by contradiction.
Assume that (\ref{qj}) holds but that
\begin{equation}
\EBjp1(t) <  q^{j+1}
\label{false}
\end{equation}
for all $t$ $\in$ $[t_0,t_0+\rho^j]$ and for all $\rho$ with $0 < \rho <1$.

Since $\EBj(t)$ is increasing with $\EBj(t_0) \ge q^j$, we have
 $$\EBj(t)=E_j(t) + \EBjp1(t) \ge q^j, \quad \forall t \ge t_0 $$
and the assumption (\ref{false}) then implies that 
\begin{equation}
E_j(t) > (1-q) q^j, \quad \forall t \in [t_0, t_0 + \rho^j].
\label{Ejbd}
\end{equation}
Now, integrating (\ref{Eflux}) written for $\EBjp1(t)$, from $t_0$ to
$t_0+\rho^j$ and using (\ref{Ejbd}) yields
\begin{equation}
\begin{split}
\EBjp1(t_0+\rho^j) &= \EBjp1(t_0)+ 2 \lambda^{j+1}
\int_{t_0}^{t_0+\rho^j} E_j a_{j+1} dt \\
&\ge  2 \lambda^{j+1} (1-q)q^j 
\int_{t_0}^{t_0+\rho^j} a_{j+1} dt.
\end{split}
\end{equation}
The assumed upper bound (\ref{false}) on $\EBjp1(t_0+\rho^j)$ then
gives
\begin{equation}
\int_{t_0}^{t_0+\rho^j} a_{j+1} dt \; \le \frac{q}{2(1-q)
  \lambda^{j+1}}.
\label{intbd}
\end{equation}
Next, consider the equation for $a_{j+1}$. Since, $E_j=a_j^2$, it reads
\begin{equation}
\frac{d a_{j+1}}{d t} 
= \lambda^{j+1} E_j  - \lambda^{j+2} a_{j+1} a_{j+2}.
\end{equation}
Integrating from $t_0$ to $t_0+\rho^j$ and using the lower bound
(\ref{Ejbd}) on $E_j$, the upper bound (\ref{intbd}) together with
the upper bound $a_{j+2} \le \sqrt{\EBjp1} < q^{(j+1)/2}$
which follows from (\ref{false}), we obtain
\begin{equation}
a_{j+1}(t_0+\rho^j) \ge q^{(j+1)/2} 
\left[\lambda^{j+1} \rho^j q^{j/2} \frac{1-q}{q^{1/2}} 
- \lambda \frac{q}{2(1-q)} \right].
\end{equation}
This is in contradiction with the assumed bound $a_{j+1}\le 
\sqrt{\EBjp1} < q^{(j+1)/2}$ provided 
\begin{equation}
\lambda \rho q^{1/2} > 1,
\label{req}
\end{equation}
for some $q$ and $\rho$ with $0<q<1$ and $0<\rho <1$,
and sufficiently large $j$. Since $\lambda >1$, it is 
always possible to pick $q$ with $\lambda^{-2} <q<1$ such that
$\lambda q^{1/2} >1 $ and then $\rho$ such that
 $\lambda^{-1} q^{-1/2} < \rho <1 $, so that condition
(\ref{req}) can always be satisfied.

\end{proof} 

Define the $H^{\alpha}$ norm of the solution as 
\begin{equation}
||a||^2_{_{H^{\alpha}}}= 
\sum_{j=j_0}^{\infty} \left(1 + (\mu^{j})^{2 \alpha}\right) |a_j|^2,
\label{Heps}
\end{equation}
where $\mu^j$ is the wavenumber associated with amplitude $a_j$,
with $\mu >1$.
For model (\ref{cmod}), we should pick $\mu=\lambda$, however 
Friedlander and Pavlovic picked $\mu=2$.
A proof of finite-time blow-up in $H^{\alpha}$ 
follows from repeated use of lemma \ref{lemma1}.

\begin{thm}
\label{thm1}
If $a_j(t)$ is a solution to system (\ref{cmod}) with $a_j(0)\ge 0$,
$\forall  j$ and 
\begin{equation}
E_{{\mathcal B}(J)}(0) \ge q^J
\end{equation}
for some sufficiently large integer $J \ge j_0$ and $\lambda^{-2}< q < 1$, 
then the $H^{\alpha}$ norm of the solution becomes unbounded
 in finite-time, for any $\alpha >0$.
\end{thm}

\begin{proof}
Applying lemma \ref{lemma1},  there is a $\rho$ with $0<\rho<1$ and
a time $t_1$ $\in$ $[0, \rho^J]$ such that 
$$E_{{\mathcal B}(J+1)}(t_1) \ge q^{J+1}.$$
 Iterating  this argument, there is a time
 $t_k$ $\in$ $[t_{k-1},t_{k-1}+\rho^{J+k-1}]$
such that $$E_{{\mathcal B}(J+k)}(t_k) \ge q^{J+k}.$$
Then 
\begin{equation}
\begin{split}
||a||^2_{H^{\alpha}} & \ge \sum_{l\ge J+k}^{\infty} \mu^{2\alpha l} a_l^2(t_k)\\
& \ge \mu^{2\alpha (J+k)} \sum_{l\ge J+k}^{\infty} a_l^2(t_k)\\
& \ge \mu^{2\alpha (J+k)}  E_{{\mathcal B}(J+k)}(t_k) \\
& \ge \mu^{2\alpha (J+k)} q^{J+k}\\
\end{split}
\end{equation}
which blows up as $k\to \infty$ provided 
\begin{equation}
\mu^{2\alpha} q > 1.
\label{bucond}
\end{equation}
This will occur in finite time since 
\begin{equation}
t_k \le \rho^J \, \left[1 + \rho + \cdots + \rho^{k-1} \right] = 
\rho^J \frac{1-\rho^k}{1-\rho},
\end{equation}
which is finite as $k \to \infty $ since $0 < \rho < 1$.

To verify that (\ref{bucond}) is compatible with the
condition (\ref{req}) required by lemma \ref{lemma1}, pick
 $q=\mu^{-2\alpha + \delta}$ with $0< \delta<2 \alpha$
so that $q<1$ and (\ref{bucond}) is satisfied, since $\mu>1$.
Condition (\ref{req}) requires
$\lambda^2 q = \lambda^2 \mu^{- 2 \alpha + \delta}  > 1
$ 
since we need $\rho^2<1$.
Hence, finite time blow-up for $\mu = \lambda$ requires
$q=\mu^{-2\alpha + \delta}$ with
\begin{equation}
\max\{0,2\alpha-2\} < \delta < 2\alpha,
\end{equation}
while if $\mu=2$ with $\lambda=2^r$, $r>0$, as in \cite{FP04}, blow-up
requires 
\begin{equation}
\max\{0,2\alpha-2r\} < \delta < 2\alpha.
\end{equation}
In either case, finite time blow-up is achievable for any $\alpha >0$.

\end{proof}

This proves that energy-conserving systems of the form (\ref{cmod}), 
with $\lambda>1$, blow-up in finite time in the $H^{\alpha}$ norm 
for any $\alpha >0$. This applies to the Friedlander-Pavlovic model,
for which $\lambda=2^{5/2}$, and is much stronger than the finite-time
blow-up in  $H^{3/2+\epsilon}$  proved in \cite{FP04}. 
Friedlander and Pavlovic concluded that if $\lambda =2^r$, then
$r>3/2$ was sufficient for blow-up, but finite time blow-up in 
$H^{\alpha}$ occurs in fact for any $r>0$. These results also apply to
the Katz-Pavlovic dyadic wavelet model which is identical, 
for special classes of  initial conditions, to the dyadic 
inviscid Burgers model (\ref{dyburg}) for which $\lambda=2$. 
Indeed, in the notation of \cite{KP04} and section \ref{KPredux}, 
with  ${\mathcal D_j}$ denoting the set of all dyadic cubes of
sidelength $2^{-j}$,
\begin{equation}
\begin{split}
||u||^2_{H^{\alpha}} & \ge 
\sum_{j=j_0}^{\infty} 2^{2\alpha j} \sum_{Q \in {\mathcal D_j}}
u^2_Q \\
& \ge \sum_{j=j_0}^{\infty} 2^{2\alpha j} \,2^{3 (j-j_0)} u^2_j \\
& =  2^{-3j_0} \sum_{j=j_0}^{\infty} 2^{2\alpha j} a_j^2
\end{split}
\end{equation}
and this last sum can diverge in finite time by theorem
\ref{thm1}. Friedlander and Pavlovic choose $\lambda = 2^{5/2}$, 
$q=2^{-3-\epsilon}$ and $\rho=2^{-\epsilon}$  with
$\epsilon >0$, so $q<1$ and $\rho<1$ as required, then 
condition (\ref{req}) reads
\begin{equation}
\lambda^2 \rho^2 q = 2^{5- 2 \epsilon - 3 - \epsilon} > 1
\end{equation} 
which requires $0 < \epsilon < 2/3$. Katz and Pavlovic, in effect,
choose $\lambda=2$, $q=2^{-\epsilon}$ and $\rho=2^{-\epsilon}$,
since they work with $u_j$ instead of $a_j=2^{3j/2} u_j$,
leading to the same restriction $0< \epsilon <2/3$.

\section{Concluding remarks}

The dyadic models of the three-dimensional Euler equations 
that have been proposed and studied by Katz and Pavlovic \cite{KP04}
and Friedlander and Pavlovic \cite{FP04} are in fact much
more closely related to the one-dimensional Burgers equation.
For the class of initial conditions discussed in section
\ref{KPredux}, the Katz-Pavlovic model is identical to a dyadic
inviscid Burgers equation. This is the model obtained by restricting
nonlinear interactions in the one-dimensional Burgers equation
to dyadic wavenumbers $k=2^j k_0$, $j$ integer. The
Friedlander-Pavlovic model reduces to a similar system but with even
stronger nonlinear couplings. The dyadic inviscid Burgers equation,
and related systems, exhibit finite time blow-up in $H^{\alpha}$,
for any $\alpha >0$. This is stronger than the finite time
blow-up for $\alpha \ge 1/2$
 that can occur in the unrestricted inviscid Burgers
equation. Nonlinear steepening is therefore enhanced by the
restriction to dyadic wavenumbers. 

The idea of studying smoothness and finite time blow-up questions
on simplified models of the Euler and Navier-Stokes equations has
merits. However, one should probably consider some of the simplified
models that have been studied in the turbulence literature 
(see \eg \cite{Bif02}). 

An early model, proposed by Obukhov \cite{O71}, bears a strong
resemblance to the dyadic model (\ref{cmod}) discussed in this paper.
The Obukhov model is
\begin{equation}
\frac{1}{\lambda^{1/3}} \frac{d a_j}{d t} 
= \lambda^{j}  a_{j-1} a_j - \lambda^{j+1}  a^2_{j+1}  - \nu_j a_j,
\label{obu}
\end{equation}
where $\nu_j \ge 0$ is a wavenumber-dependent viscosity coefficient,
for instance $\nu_j= \nu \lambda^{2j}$, with $\nu >0$ to model Navier-Stokes. 
The dynamically inconsequential $\lambda^{-1/3}$ factor
on the left-hand side has been introduced so that both models
 (\ref{cmod}) and (\ref{obu}) have the same power-law steady states 
\begin{equation}
a_j = \lambda^{-2/9} {{\mathcal E}}^{1/3}\, \lambda^{-j/3}
\label{power}
\end{equation}
in the unbounded, inviscid limit ($\nu=0$, $j_0 \to -\infty$), where
 \begin{equation}
{\mathcal E}= \frac{1}{2} \frac{d}{dt} \EBj,
\end{equation}
is the energy flux which is independent of $j$ for the power-law (\ref{power}).
This power law is analogous to the Kolmogorov-Obukhov $k^{-5/3}$
power law for the energy spectrum in three-dimensional Navier-Stokes
turbulence.
However, the dynamics of the two models are quite different. Model (\ref{cmod}) is a
nonlinear steepening model where energy is pushed ever more
efficiently to larger wavenumbers. Model (\ref{obu}) is an instability
cascade model. Consider, for example, initial
conditions such that 
 all modes are initially zero except $a_{l}(0)>0$,
for some integer $l$. The nonlinear steepening model (\ref{cmod})
 directly starts pushing energy to $j>l$, but the energy would cascade to $j<l$
in the Obukhov model. 
The energy will cascade to $j=l+1$ only if 
$a_{l+1}(0) \ne 0$ and $\lambda^{l+1} a_l > \nu_{l+1}$.
The cascade will continue if $a_{l+2}$ is non-zero and $\lambda^{l+2} a_{l+1} >
\nu_{l+2}$, or in other word if the local Reynolds number is larger
than 1. It would be interesting to show what type of finite time
blow-up is possible in the inviscid Obukhov model and see whether
there is any connection between blow-up and the power-law
(\ref{power}). 
A connection is expected because $a_j>0$ is required for energy cascade to  
larger $j$ and $a_j(t)$ is guaranteed to stay positive as long as
$a_{j-1}a_j \ge \lambda a_{j+1}^2$, when $\nu=0$. For power law scaling, 
$a_j \propto \lambda^{\beta j}$, this requires $\beta \le -1/3$.

\end{document}